# Apex control within an elasto-plastic constitutive model for confined concretes

G. Mazzucco[1*], B. Pomaro[1], V.A. Salomoni[1,2], C.E. Majorana[1]

[1] Department of Civil, Environmental and Architectural Engineering (DICEA), University of Padua, Via F. Marzolo 9, 35131 Padova, Italy

[2] Department of Technology and Management (DTG), University of Padua, Stradella S. Nicola 3, 36100 Vicenza, Italy

**Abstract**

This work focuses on the numerical modelling of confined concretes when subjected to multiaxial loading accounting for a nonlinear material response. Particularly, an improved return mapping procedure is here presented to cope with the issue of singularities (corners) in the yield surface, able to catch locally the optimal return point on the active yield surface. The algorithm is applied to the Menétry and Willam's yield surface, with a plastic potential as suggested by Grassl and the softening law proposed by Cervenka suitable for concrete materials. The model capabilities to reproduce the confined behaviour of concrete, as well as its post-peak behaviour, and to predict perfect plasticity or softening are here discussed. The proposed approach can be applied to any elastoplastic material with singular yield surface along the hydrostatic axis.

**Keywords:** computational plasticity; multiaxial stress state; softening; apex; yield surface.

[*] Corresponding author.
Email gianluca.mazzucco@unipd.it
Ph. +39 049 8275565

**Nomenclature**

| | |
|---|---|
| *a* | *Apex point* |
| *B* | *Modified bulk modulus* |
| $B^*$ | *Bulk modulus* |
| **C** | *Elastic tangent operator* |
| $C^{ep}$ | *Elastoplastic tangent operator* |
| *E* | *Young modulus* |
| *e* | *Material eccentricity* |
| *f* | *Yield surface* |
| $f_c$ | *Compression strength* |
| $f_t$ | *Tensile strength* |
| *g* | *Plastic potential* |
| $g_A$, $g_B$ | *Constant parameters in the plastic potential* |
| *κ* | *Equivalent plastic strain (volumetric plastic strain)* |
| $k_{1D}$ | *Plastic volumetric strain at maximum compression load in uniaxial condition* |
| **I** | *Identity vector $[1,1,1,0,0,0]^T$* |
| $I_1$ | *First stress tensor invariant* |
| $J_2$, $J_3$ | *Second and third deviatoric invariants* |
| *q* | *Hardening/softening parameter; $q=q_h \cdot q_s$* |
| $q_h$ | *Hardening parameter* |
| $q_{h0}$ | *Initial hardening parameter* |
| $q_s$ | *Softening parameter* |
| *r* | *Error or residual variable* |
| *s* | *Deviatoric stress tensor* |
| *t* | *Softening slope parameter* |
| $\varepsilon^e$ | *Elastic strain tensor* |
| $\varepsilon^p$ | *Plastic strain tensor* |
| $\varepsilon_{ij}$ | *Total strain component* |
| $\varepsilon^e_{ij}$ | *Elastic strain component* |
| $\varepsilon^p_{ij}$ | *Plastic strain component* |
| *δ* | *Kronecker delta ($\delta_{ij}=1$ if $i=j$; $\delta_{ij}=0$ if $i \neq j$)* |
| *σ* | *Stress tensor* |

| $\bar{\sigma}$ | *Hydrostatic stress* |
| --- | --- |
| *v* | *Poisson's coefficient* |
| *ξ$_a$* | *Coordinate ξ in the apex point* |
| *ξ$^{Tr}$* | *Coordinate ξ for the trial stress* |
| *(ξ, ρ, θ)* | *Haigh-Westergaard coordinates* |

**Introduction**

Several studies have addressed the role of confinement in concrete, leading to formulations for the behaviour of concrete in multiaxial stress states [1]-[6], sometimes including damage and fracture [7]-[11].

As a path-dependent material, the updating scheme requires the formulation of a numerical algorithm for integration of the corresponding rate constitutive equations. There are several methods to solve such an issue [12]; typically, an Euler difference scheme is applied to find an approximate solution to the elastoplastic problem, in the form of a two-step algorithm: the *elastic predictor* (in which the evolution problem is solved as if the material were purely elastic in the current interval), followed by the *return mapping or plastic corrector* [13] (which accounts for plastic flow and enforces plastic admissibility).

In other words, the problem consists in determining the stress on the failure surface satisfying the yield condition as well as the hardening/softening law.

The scheme is relatively straightforward for a smooth yield surfaces; unfortunately, in tensile conditions elastoplastic models for concretes (even not confined) might present some limitations coming from the use of a non-associated flow rule. In such cases, where the form of the plastic potential differs from that of the yield surface (which implies that the plastic flow is not perpendicular to the yield surface), the return mapping procedure does not guarantee a solution on the yield surface due to the shape of the plastic potential, namely the presence of singularities or *apex* points.

To overcome these restrictions an improved return mapping algorithm has been implemented to solve the local problem connected to the singularity of the Menétry and Willam's yield surface [1] under tensile stress states, in line with the theory of plasticity [12],[14]. This surface, formerly proposed by Etse and Willam [15], is established in terms of cylindrical coordinates in the principal effective stress space (Haigh-Westergaard coordinates) and it is one of the most suitable smooth surfaces to describe the triaxial failure envelope of concrete, covering in a unified way the tension-tension, tension-compression and compression-compression mechanisms. The model here studied is completed by the plastic potential suggested by Grassl and co-authors [4], which is based on the adoption of the volumetric plastic strain as hardening parameter.

In the first Section of the paper, the theoretical aspects related to the upgraded return-mapping procedure are provided, including the numerical strategies necessary to enforce the return to the apex of the trial stress. At the same time a check for activating the appropriate return mapping method is included. The subsequent Section presents some benchmark tests for validating the algorithm under different tensile stress states.

## 1 Theoretical formulation

The yield function $f$ and the plastic potential $g$ are defined through the Haigh-Westergaard coordinates ($\xi, \rho, \theta$):

$$\xi = \frac{1}{\sqrt{3}} I_1 \; ; \; \rho = \sqrt{2 J_2} \; ; \; \cos 3\theta = \frac{3\sqrt{3}}{2} \frac{J_3}{J_2^{3/2}} \tag{1}$$

hence, referring to Menétrey and Willam [1]:

$$f = \frac{3}{2}\left(\frac{\rho}{f_c}\right)^2 + q_h m \left(\frac{\rho}{\sqrt{6} f_c} R + \frac{\xi}{\sqrt{3} f_c}\right) - q \leq 0. \tag{2}$$

where $m$ is the cohesion parameter, depending on the uniaxial strength of concrete in compression and tension $f_c$ and $f_t$, respectively, through the eccentricity $e$, defined as of out-of-roundness of the yield surface in the deviatoric section:

$$m = \frac{3\left(f_c^2 - f_t^2\right)}{f_c f_t} \frac{e}{e+1} \tag{3}$$

while $R = R(\theta, e)$ is the elliptic function, which describes the roundness of the failure surface (see [4],[9]):

$$R = \frac{4\left(1-e^2\right)\cos^2\theta + \left(2e-1\right)^2}{2\left(1-e^2\right)\cos\theta + \left(2e-1\right)\sqrt{4\left(1-e^2\right)\cos^2\theta + 5e^2 - 4e}} \tag{4}$$

The hardening/softening parameter $q$ is dependent on the equivalent strain $\kappa$ and it can be split into two components: the hardening part $q_h$ and the softening part $q_s$:

$$q(\kappa) = q_h(\kappa) \, q_s(\kappa). \tag{5}$$

The Menétrey-Willam failure surface is not fixed but it can move along the hydrostatic axis based on the value of $\kappa$. Coherently with the novel hardening law reported in [4] this factor coincides with the volumetric plastic strain $\varepsilon_v^p$ with the advantage that the material behaviour in a multiaxial regime can be calibrated via a reduced number of parameters.

The formulation of the $q$ function for concrete materials has been suggested by Červenka and Papanikolaou [9]:

$$q_h(\kappa) = q_{h0} + (1 - q_{h0})\sqrt{1 - \left(\frac{k_{1D} - \kappa}{k_{1D}}\right)^2}$$

$$q_s(\kappa) = \left[1 + \left(\frac{n_1 - 1}{n_2 - 1}\right)^2\right]^{-2}$$

(6)

where $q_{h0}$ is the value of $q_h$ when $\kappa=0$, which defines the yield surface at the initial elastic regime; $k_{1D}$ is the plastic volumetric strain at uniaxial concrete strength; $n_1=\kappa/k_{1D}$; $n_2=(k_{1D}+t)/k_{1D}$; $t$ is the slope of the softening function. A graphical representation of $q_h$ and $q_s$ is reported in Fig. 1. It can be seen that the hardening function is constantly equal to unit at the end of the hardening process, when the material undergoes softening; the softening parameter starts from unit and decreases to zero at increasing decohesion of the material. At the same time the failure surface shifts along the negative hydrostatic axis.

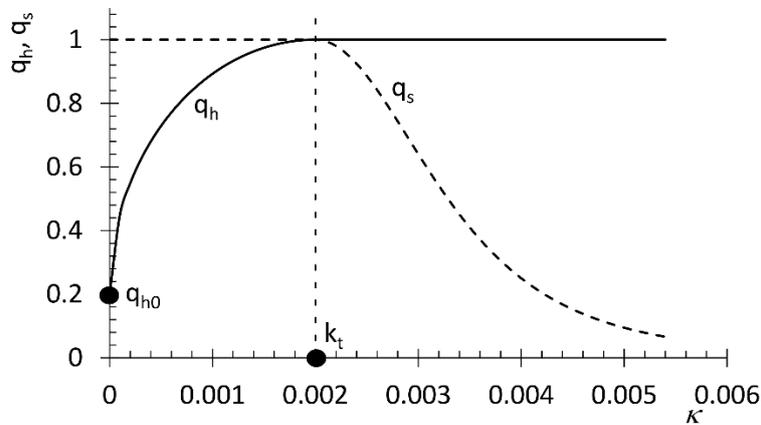

*Fig. 1 Hardening/softening parameters vs. $\kappa$*

The adopted plastic model is non-associated with the plastic potential being the quadratic expression proposed in [4]:

$$g = -g_A\left(\frac{\rho}{f_c\sqrt{q}}\right)^2 - g_B\frac{\rho}{f_c\sqrt{q}} + \frac{\xi}{f_c\sqrt{q}} = 0 \tag{7}$$

and the two parameters $g_A$ and $g_B$ can be determined starting from the axial strain at maximum stress in uniaxial compression and in triaxial compression tests.

For the integration of the constitutive equations, an implicit backward-Euler return mapping algorithm is used. This method implies that the return mapping vector is obtained by evaluating the normal vector to the plastic potential $g$ computed in the trial point.

Anyway, if the described procedure is correctly applied to concrete behaviour in compression, in traction the method could fail [4] since the normal vector to the $g$ function $\partial g / \partial \boldsymbol{\sigma}$ cannot generally intercept the yield function for a trial stress $\sigma^{Tr}$ overcoming the apex point $a$ (see Fig. 2).

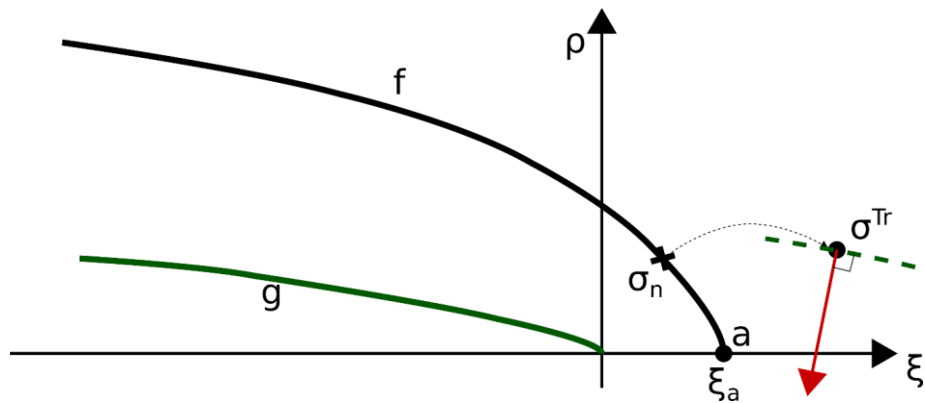

*Fig. 2 Classic return mapping procedure in a traction condition.*

In fact under such circumstances convergence does not occur unless an enhanced return mapping procedure, which allows the return to the apex of the yield surface, is adopted.

Therefore, a "return-to-apex procedure" has been innovatively introduced in the numerical integration scheme, in line with the guidelines reported in [12].

According to this procedure, the return mapping is conceived as a projection of the trial stress onto the hydrostatic direction with the aim to have the final stress $\sigma_{n+1}$ equal to the point $a$ ($\xi_{n+1}=\xi_a$, see Fig. 3). It can be noticed that the apex point is a function of the equivalent plastic volumetric strain $\kappa$: $a=a(\kappa)$, so a nonlinear procedure must be accomplished to get the solution. Furthermore, the technique to check whether to activate or not the return-to-apex control is straightforward: it consists on controlling whether $\xi^{Tr}$ overcomes or not the limit in traction of the yield surface $f$, i.e. $\xi_a$ (see Fig. 3).

More in detail, the procedure is outlined in the following. Point $a$, which is the point of maximum hydrostatic tensile stress, is characterized by two conditions on the yield function: $f=0$ and $\rho=0$, therefore:

$$\xi_a(\kappa) = \frac{\sqrt{3} f_c}{m} q_s(\kappa). \tag{8}$$

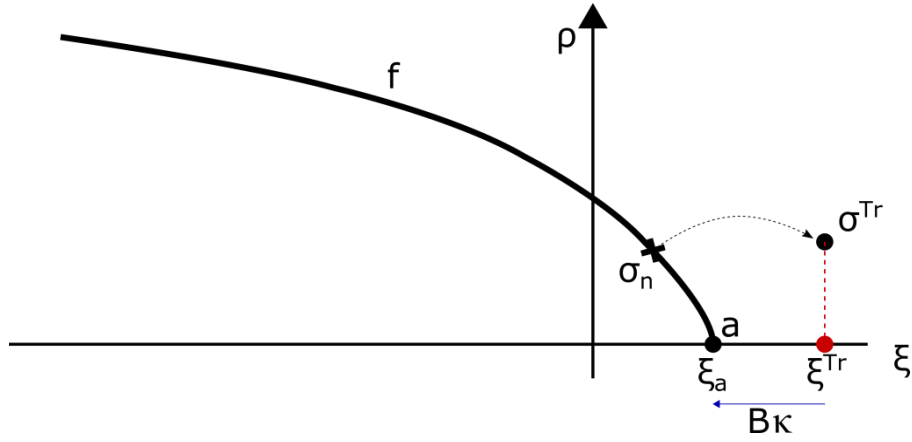

*Fig. 3 Modified return mapping procedure.*

If $\xi^{Tr} > \xi^a$ then at the apex the return mapping procedure is given by Perić and Owen [12] along the hydrostatic axis:

$$\xi_{n+1} = \xi^{Tr} - B\kappa \tag{9}$$

where $B$ is a modified bulk modulus. In fact, considering that $\xi = I_1/\sqrt{3}$ and that the bulk modulus, defined by Eq. (10a), relates the hydrostatic stress with the volumetric strain (10b):

$$B^* = \frac{E}{3(1-2\nu)} \quad a) \quad \text{and} \quad \kappa = \frac{I_1}{3B^*} \quad b) \tag{10}$$

where $I_1/3 = \bar{\sigma}$ is the hydrostatic stress, then the modified bulk modulus $B$ stands in the following relationship with the bulk modulus $B^*$:

$$\kappa = \frac{\sqrt{3}\xi}{3B^*} = \frac{\xi}{B} \quad \text{with} \quad B = \sqrt{3}B^* = \frac{E}{\sqrt{3}(1-2\nu)}. \tag{11}$$

Since $q_s$ is function of $\kappa$, the solution of Eq. (9) is in general affected by an error on the residual $r$:

$$r = \xi_{n+1} - \left[\xi^{Tr} - B\kappa\right]. \tag{12}$$

But when a solution $\xi_{n+1}$ is found coincident with the apex, this means that $\xi_{n+1} = \xi_a$ and Eq. (12) can be rewritten as:

$$r(\kappa) = \xi_a(\kappa) - \left[\xi^{Tr} - B\kappa\right]. \tag{13}$$

To minimize the error $r$ an implicit backward-Euler scheme is applied on the error variation $\partial r$:

$$\partial r = \partial \xi_a - \left(\partial \xi^{Tr} - B\partial\kappa\right) = \frac{\sqrt{3}f_c}{m}\frac{\partial q_s}{\partial \kappa} + B. \tag{14}$$

The return mapping algorithm enhanced by the return-to-apex procedure is summarized in Box. 1.

**Box 1.**

The return-to-apex procedure at a step *n+1* is:

$\xi^{Tr} > \xi_a(\kappa_n)$

Variable initialization:

$^0\kappa = \kappa_n$ (where $\kappa_n$ is the equivalent plastic strain at step *n*)

$^0r = r\left(^0\kappa\right) = \xi_a - \left(\xi^{Tr} - B{}^0\kappa\right)$

**for** *m=1; MAXITER*

$\quad ^m\kappa = {}^{m-1}\kappa + \dfrac{^{m-1}r}{\partial^{m-1}r}$

$\quad q_s\left(^m\kappa\right)$ (see eq. 6); $\quad ^m\xi_a = \dfrac{\sqrt{3}f_c}{m}q_s\left(^m\kappa\right)$

$\quad ^mr = r\left(^m\kappa\right) = \xi_a - \left(\xi^{Tr} - B{}^m\kappa\right)$

$\quad \partial^m r$ (see eq. 14)

$\quad$ **if** *r < TOLL;*

$\quad\quad \kappa_{n+1} = {}^m\kappa; \quad \xi_{n+1} = {}^m\xi_a$

$\quad\quad$ **exit**

$\quad$ **end if**

**end for**

When the convergence is obtained (at step *n+1*) the stress tensor results in a hydrostatic configuration with:

$$\bar{\sigma}_{n+1} = \frac{\xi_{n+1}}{\sqrt{3}} \tag{15}$$

and the stress tensor **σ** can be obtained as $\boldsymbol{\sigma}_{n+1} = \bar{\sigma}_{n+1}\mathbf{I}$ with **I** identity tensor.

The updated plastic strain tensor $\boldsymbol{\varepsilon}^p_{n+1}$ can be found starting from the updated stress tensor:

$$\boldsymbol{\varepsilon}_{n+1}^{p} = \boldsymbol{\varepsilon}_{n}^{p} + \mathbf{C}^{-1}\left(\boldsymbol{\sigma}^{Tr} - \boldsymbol{\sigma}_{n+1}\right). \tag{16}$$

It is to be noticed that at the apex the constitutive tangent operator must be recalculated ($C^{ep}$) taking into account the final configuration of the stress tensor $\boldsymbol{\sigma}_{n+1}$:

$$\mathbf{C}^{ep} = \frac{\partial \boldsymbol{\sigma}_{n+1}}{\partial \boldsymbol{\varepsilon}_{n+1}^{e,Tr}} = \frac{1}{\sqrt{3}} \frac{\partial \xi_{n+1}}{\partial \boldsymbol{\varepsilon}_{n+1}^{e,Tr}} \mathbf{I} \tag{17}$$

where $\mathbf{I}\boldsymbol{\varepsilon}_{n+1}^{e,Tr}$ is the trial strain tensor at step $n+1$, i.e. the elastic component at the same step; while $\partial \xi_{n+1}$ is computed as:

$$\partial \xi_{n+1} = \partial \xi^{Tr} - B \partial \kappa_{n+1} = B \, \mathbf{I} \partial \boldsymbol{\varepsilon}_{n+1}^{e,Tr} - B \partial \kappa_{n+1} \tag{18}$$

in which $\mathbf{I} \partial \boldsymbol{\varepsilon}_{n+1}^{e,Tr} = \kappa_{n+1}^{e,Tr}$ is the volumetric elastic strain in the trial state. In order to obtain $\partial \kappa_{n+1}$ Eq. (14) must vanish, i.e. $\partial \kappa_{n+1}$ is found by imposing $\partial r = 0$, which occurs when the solution is found:

$$\partial r = \frac{\sqrt{3} f_c}{m} \frac{\partial q_s}{\partial \kappa} \partial \kappa_{n+1} - B \, \mathbf{I} \partial \boldsymbol{\varepsilon}_{n+1}^{e,Tr} + B \partial \kappa_{n+1} = \left(\frac{\sqrt{3} f_c}{m} \frac{\partial q_s}{\partial \kappa} + B\right) \partial \kappa_{n+1} - B \, \mathbf{I} \partial \boldsymbol{\varepsilon}_{n+1}^{e,Tr} = 0 \tag{19}$$

that yields:

$$\partial \kappa_{n+1} = \frac{B}{\frac{\sqrt{3} f_c}{m} \frac{\partial q_s}{\partial \kappa} + B} \mathbf{I} \partial \boldsymbol{\varepsilon}_{n+1}^{e,Tr}. \tag{20}$$

At this point Eq. (17) can be rewritten as:

$$\partial \xi_{n+1} = B \, \mathbf{I} \partial \boldsymbol{\varepsilon}_{n+1}^{e,Tr} - \frac{B^2}{\frac{\sqrt{3} f_c}{m} \frac{\partial q_s}{\partial \kappa} + B} \mathbf{I} \partial \boldsymbol{\varepsilon}_{n+1}^{e,Tr} \tag{21}$$

and the consistent tangent operator at the apex can be computed as:

$$\mathbf{C}^{ep} = \frac{B}{\sqrt{3}} \left[ 1 - \frac{B}{\left(\frac{\sqrt{3} f_c}{m} \frac{\partial q_s}{\partial \kappa}\right) + B} \right] \mathbf{I} \otimes \mathbf{I}. \tag{22}$$

## 2 Model validation

Some examples have been developed to check the correctness of the procedure and verify the response of the model under different dominant tensile loading histories. For this purpose on the same reference sample the material characteristics or the load history have been slightly modified, so to appreciate the nonlinear mechanical response of the material under different scenarios.

## 2.1 Hydrostatic tensile stress state with constant $q_s$

A cubic sample with side length of *100mm* has been loaded under a hydrostatic traction condition. In this example the softening parameter $q_s$ has been maintained constantly equal to *1* (perfect plasticity). This condition is reached by increasing $k_{1D}$ considerably more than the volumetric strain obtained during the whole analysis.

The material parameters and the constitutive characteristics of the model have been reported in Tab. 1. Let's observe that, with $q_s$ equal to unit during the whole analysis, the material behaviour is independent on the softening parameter slope *t*.

*Tab. 1 Material and constitutive coefficients.*

| | | |
|---|---|---|
| E | 30000.00 | MPa |
| v | 0.15 | |
| fc | 32.00 | MPa |
| ft | 3.00 | MPa |
| e | 0.52 | |
| t | 0.0055 | |
| $k_{1D}$ | 0.10008 | |
| $q_{h0}$ | 0.20 | |
| $g_A$ | 21.22 | |
| $g_B$ | 31.46 | |

By assuming that $q_s$ is constant during the stress evolution means that the apex point remains constant (no shifting along the hydrostatic axis) together with $\xi_a$, (Fig. 4). At increasing stresses, the hardening parameter $q_h$ increases with a perfect plasticity behaviour, softening not occurring. The elastic limit in traction is $\xi_a$=5.11MPa, which can be found from:

$$\xi_a = \frac{\sqrt{3}f_c}{m} \quad \text{with} \quad \bar{\sigma} = \frac{\xi_a}{\sqrt{3}} \tag{23}$$

and the stress tensor in the plastic condition results hydrostatic with $\bar{\sigma}$=2.95MPa. The corresponding stress vs. strain curve is shown in Fig. 5.

## 2.2 Hydrostatic tensile stress state with variable $q_s$

This example assumes the same geometry, boundary conditions, material and constitutive characteristics as those of the previous one; a reduced $k_{1D}$ only has been considered ($k_{1D}$=0.0001). Such an assumption leads to obtain a softening parameter lower than unit during the analysis, which allows the triggering and spreading of softening itself. Fig. 6a shows the evolution of $q_s$: it is equal to *1* when $\kappa \leq k_{1D}$ and it decreases when $\kappa$ overcomes $k_{1D}$ (softening). In the ($\xi, \rho$) plane, when softening occurs the apex point of the yield function eventually shifts to *0* when $q_s(\kappa)=0$, as shown in Fig. 6b.

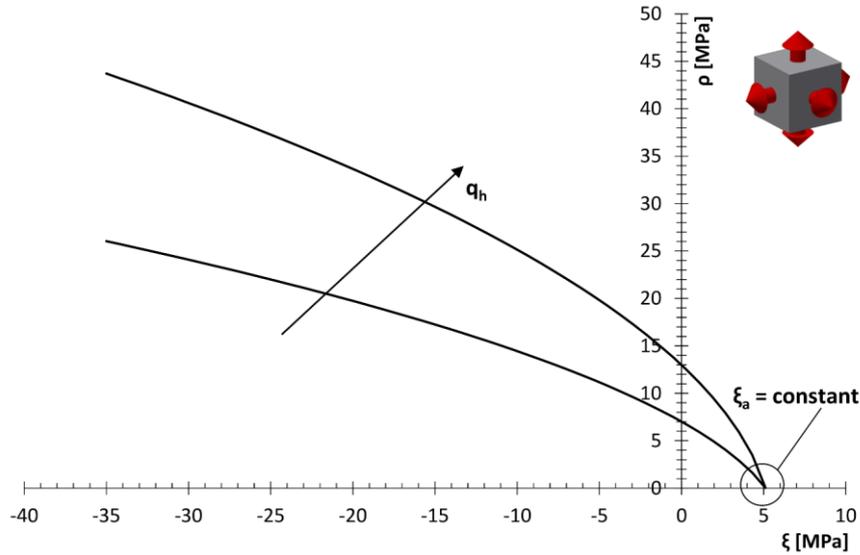

Fig. 4 Variation of the yield function in case of hydrostatic traction with $q_s$=1.

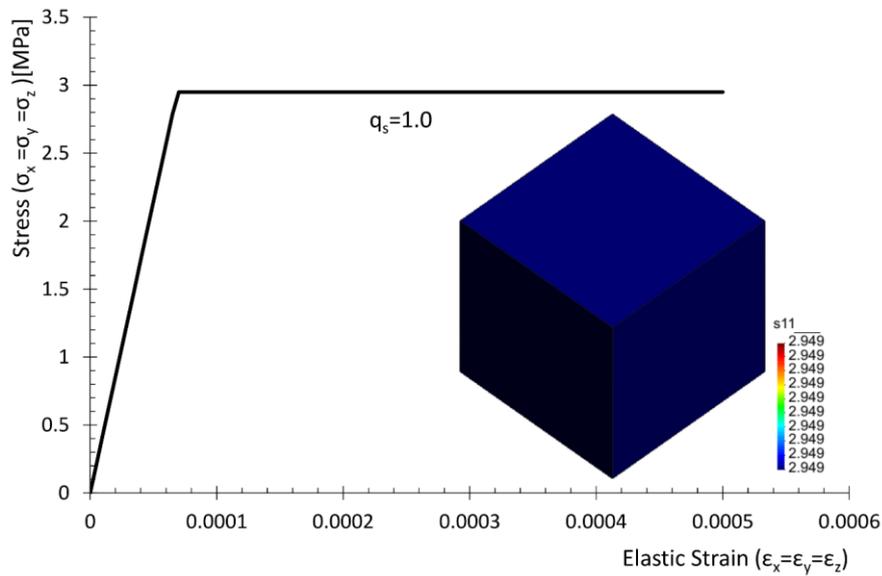

Fig. 5 Stress vs. strain in case of hydrostatic traction with $q_s$=1.

By considering the stress vs. strain curve at a reference Gauss point of the model (Fig. 7), after an elastic behaviour up to $\bar{\sigma}$ =2.95MPa, stress remains constant (condition of perfect plasticity) until $\kappa<k_{1D}$ when $q_s(\kappa)$ is still unit (up to a total strain $\varepsilon_{ii} = \varepsilon_{ii}^e + \varepsilon_{ii}^p$ equal to $k_{1D}$ =0.0001), subsequently the volumetric strain overcomes $k_{1D}$ and softening occurs.

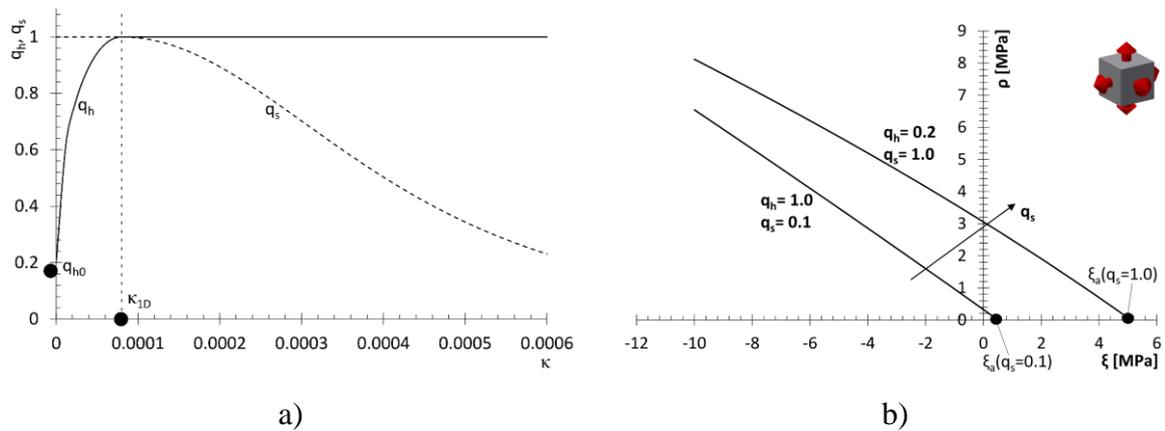

*Fig. 6 Case of hydrostatic traction with variable $q_s$ – hardening/softening parameter vs. volumetric strain a); yield surface at different $q_s$ values b).*

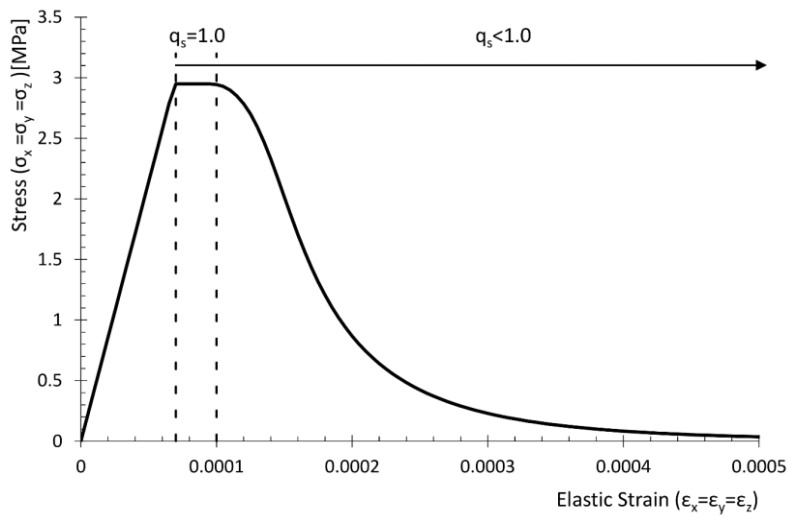

*Fig. 7 Stress vs. strain in case of hydrostatic traction with variable $q_s$.*

In case of a loading/unloading scenario, unloading typically occurs elastically (Fig. 8) [14].

## 2.3 Pre-confined stress state with constant $q_s$

In the previous examples the only possible solution was the apex one and they served just to test the correctness of the upgraded algorithm; then, more general stress states have been taken into account so to appreciate the suitability of the whole procedure.

The sample (with the same material and constitutive parameters as before; $q_s = 1$) has been hence subjected to a compressive hydrostatic stress of *-8.0MPa* ($\xi_0=13.86MPa$) (Fig. 9b) and subsequently to a variable tensile stress state (Fig. 9c); the latter scenario has been obtained by applying a variable displacement history (Fig. 9d): one along the axial direction (*z* coordinate in Fig. 9a), with maximum displacement of *0.05mm*, and one along the transversal direction *x*,

*y* (maximum displacement of *0.02mm*). In this way, as evidenced below, a deviatoric stress state triggers and spreads up to and after yielding.

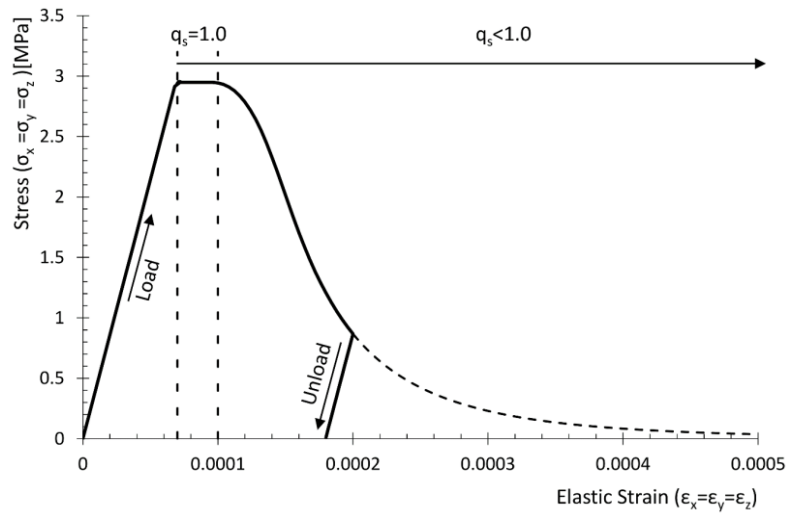

*Fig. 8 Stress vs. strain during loading/unloading in case of hydrostatic traction with variable $q_s$.*

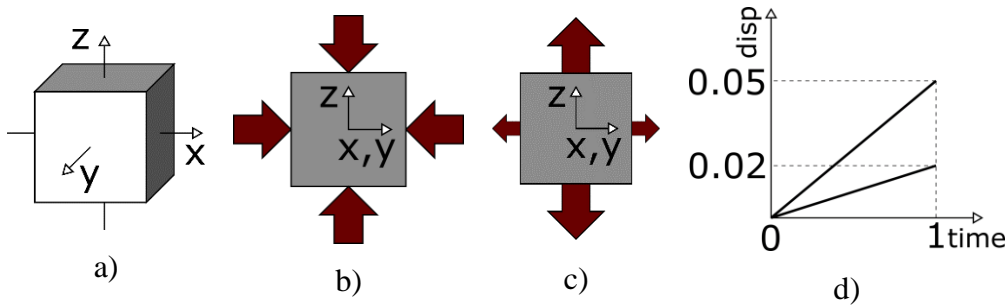

*Fig. 9 Cubic sample a); initial confinement b); deviatoric tensile state c); displacement increment d).*

Stress-strain curves along axial and transversal direction (Fig. 10) and stress path in the $(\xi, \rho)$ plane (Fig. 11) evidence the anticipated behaviour: starting from a confined state, the yield surface is reached when $\rho$ is about *2.6MPa*, then plastic deformations develop and the standard Newton-Raphson algorithm applies with decreasing $\rho$ until being close to the apex, where the return-to-apex condition is activated. When the apex is reached a tensile hydrostatic state occurs and the axial and the transversal stresses take the same values again (blue symbols in Fig. 10).

In order to check the correct activation of the return-to-apex algorithm, just the last trial stress has been plotted in Fig. 11 (red circle); this stress level is found to occur at about *50%* of the analysis and the trial configuration $\xi^{Tr}$ clearly overcomes the apex point ($\times$ symbol in Fig. 11), which definitely requires the activation of the enhanced return mapping procedure. The stress is hence returned to and maintained on the apex point.

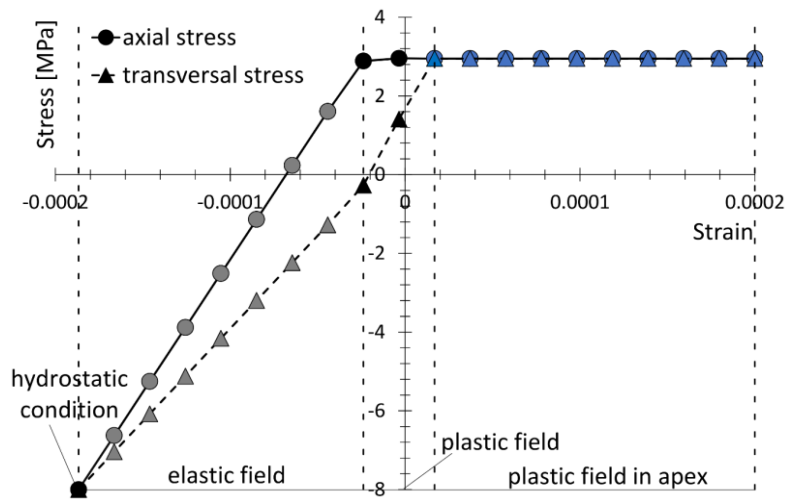

Fig. 10 Stress vs. strain along axial and transversal directions in case of pre-confinement, tensile deviatoric stress state and constant $q_s$.

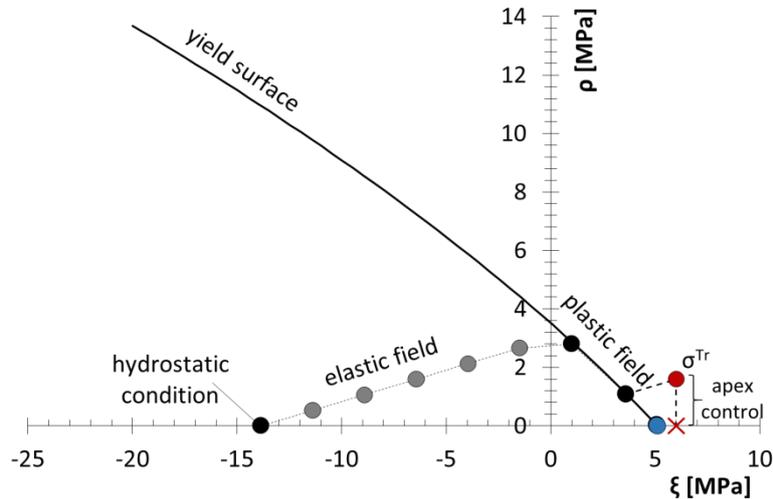

Fig. 11 Stress path in the case of pre-confinement, tensile deviatoric stress state and constant $q_s$.

### 2.4 Pre-confined stress state with variable $q_s$

In analogy with §2.2 a pre-confined state with variable softening parameter $q_s < 1$ has been considered and, as performed in §2.3, a compressive hydrostatic stress of -8.0MPa ($\xi_0=13.86MPa$) has been followed by a variable tensile stress state. The maximum axial (0.075mm)-transversal (0.005mm) displacements ratio is now increased to enhance softening, as well as $k_{1D}$ decreased ($k_{1D}=0.00008$) together with the softening parameter slope $t$ ($t=0.001$). By analysing stress-strain curves along axial and transversal direction (Fig. 11) and stress path in the ($\xi, \rho$) plane (Fig. 11), a softening behaviour is evidenced and correspondingly the apex point of the yield function eventually shifts along the negative hydrostatic axis (dotted line).

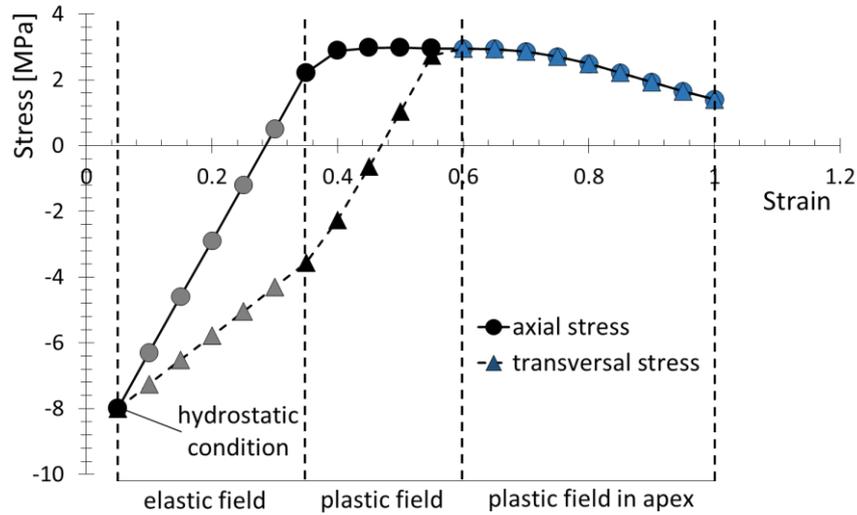

*Fig. 12 Stress vs. strain along axial and transversal directions in case of pre-confinement, tensile deviatoric stress state and variable $q_s$.*

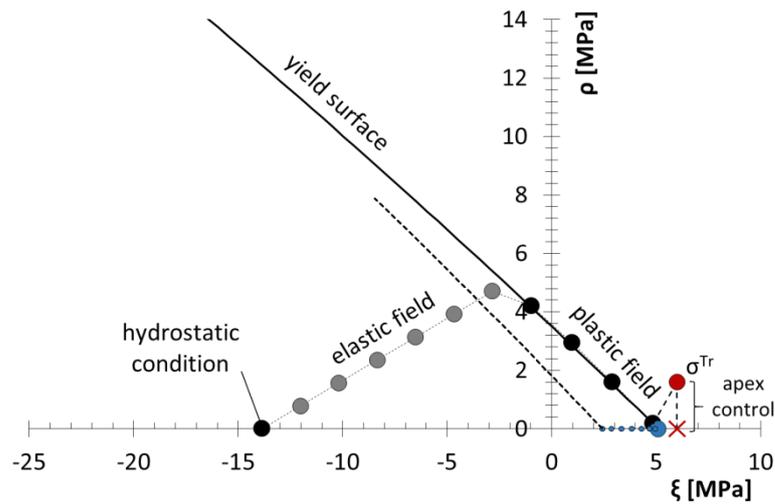

*Fig. 13 Stress path in the case of pre-confinement, tensile deviatoric stress state and variable $q_s$.*

**Conclusions**

The work presents a numerical procedure developed to overcome non-convergence issues related to tensile stress states when using a non-associated flow rule within plasticity models. This method consists in the activation of an enhanced return mapping procedure in proximity of the yield surface singularity (apex), by boosting a return-to-the apex condition. The procedure has been applied to the Menetréy and Willam yield surface with the non-associated flow rule suggested by Grassl and co-authors, particularly referring to concrete materials in confined states. The response of some reference cubic samples has been considered, allowing for both checking the correctness of the proposed algorithm and evidencing the material response under a variety of load cases, passing from hydrostatic tensile to hydrostatic compressive (confined) and deviatoric tensile stress states; a perfectly plastic or softening behaviour has been obtained. Even if here applied to specific non associated plasticity models, the procedure does not lack in generality and can be directly applied to any elastoplastic model presenting singularities along hydrostatic axes without restrictions.

**Acknowledgments**

Financial support from the Italian Ministry of Education, University and Research (MIUR) in the framework of the Project PRIN "COAN 5.50.16.01"- code 2015JW9NJT - is gratefully acknowledged.